\documentclass[12pt,fleqn]{article}
\usepackage{graphicx}

\begin{document}
\title{\textbf{A Link Between The Continuous And The Discrete Logistic Equations}}
\author{B.G.Sidharth and B.S.Lakshmi\footnote{Address for communication :International
Institute for Applicable Mathematics and Information Sciences,
B.M. Birla Science Centre, Adarsh Nagar, Hyderabad - 500
063(India)}}
\date{}
\maketitle \begin{abstract} Two types of population models are
well known -- the continuous and the discrete types.The two have
very different characteristics and methods of solutions and
analysis.In this note, we point out that an iterative technique
when applied to the continuous case mimics, surprisingly the
discrete theory. The implication is that techniques and
conclusions of the latter theory can now be applied to the former
case (and vice versa).
\end{abstract}
\section{Discrete Logistic Equation} Population growth in nature
is seldom as smoothly continuous as a classical logistic curve
suggests \cite{freedman}. In a species with a short annual
breeding season whose members live for several breeding seasons
and die at any time of the year, a continuous record of population
size would undoubtedly show seasonal undulations. For many species
in fact population growth is markedly discontinuous. These are,
for example species whose  members reproduce only once in their
lifetime and die before their descendants' lives begin, that is
each generation dies before the eggs are hatched,or the seeds are
germinated, to form the successive generations. A continuous-time
differential equation is then inappropriate as a representation of
population growth.We have to consider instead difference
equations. As this will be needed in the sequel, we touch upon the
salient features. \subsection{First-Order Difference Equation} One
of the simplest systems an ecologist can study is a seasonal
breeding population in which generations do not overlap. Many
natural populations,particularly among temperate zone insects
(including the many economically important crop and orchard pests)
are of this kind. In this situation, the observed data will
consist of information about the maximum or the average, or the
total population in each generation \cite{hallam}. The studies try
to form a relation between the magnitude of the population in the
generation \textbf{n} represented by $x_{n}$ and the magnitude of
the population in generation \textbf{n+1} represented by $x_{n+1}$
such a relation may be expressed in the general form
\begin{equation}\label{may1}
x_{n+1}= F(x_{n})
\end{equation}
The function $ F(x_{n})$ will be what the biologist calls
``density dependant", and a mathematician calls non-linear
\cite{may}. Equation (\ref{may1}) is thus a first order,
non-linear difference equation. Equation (\ref{may1}) also
describes many other examples
in biology apart from population growth. \\For instance:\vspace{4mm}\\
1. In genetics, where the equation describes the change in gene
frequency in time. \vspace{4mm}\\ 2. In epidemiology  where
$x_{n}$ represents the fraction of population infected in time n.
\vspace{4mm} \\3. In economics where the relationships between
commodity quantity and price are studied.  \vspace{4mm} \\4. In
social sciences to study the propagation of rumors where $x_{n}$
could be the number of people to have heard the rumor after time
t.  \\In many of these contexts, and for biological populations in
particular, there is a tendency for the variable $ x_{n}$ to
increase from one generation to the next when it is small, and for
it to decrease when it is large. That is, the nonlinear function
$F( x_{n}) $ often has the following properties: $ F(0)=0$ and
$F(x)$ increases monotonically as x increases through the range
$0<x< A$ with F(x) attaining its maximum value at x = A, and F(x)
decreases monotonically as x increases beyond x = A. A specific
example is afforded by the equation
\begin{equation}\label{may2}
F( x_{n})=x_{n+1}= x_{n}(a-b x_{n})
\end{equation}

This is sometimes referred to as the logistic difference equation.
Using the substitution $x_{n}=b/a  x_{n}$ Equation (\ref{may2})
can be written as
\begin{equation}\label{may3}
F( x_{n})=x_{n+1}=a x_{n}(1- x_{n})
\end{equation} The behavior of the
solutions of (\ref{may3}) is a function of the parameter `a'
Equation (\ref{may3}) has meaningful solutions for
\[ 0 \leq a\leq 4 \] with  $x_{n}$
measuring a non-negative quantity. \\Studies have shown \cite{may}
that the very simple nonlinear difference equation can possess an
extraordinarily rich spectrum of dynamical behavior, from stable
points, stable cycles to ultimately chaotic behavior.Thus the
problem is far richer than the continuous case seen earlier. In
this form equation (\ref{may3}) is a simple nonlinear difference
equation.\section{Examples of chaotic dynamical systems : The
logistic map}\label{chatd} 
\begin{figure}
\centering
 \includegraphics{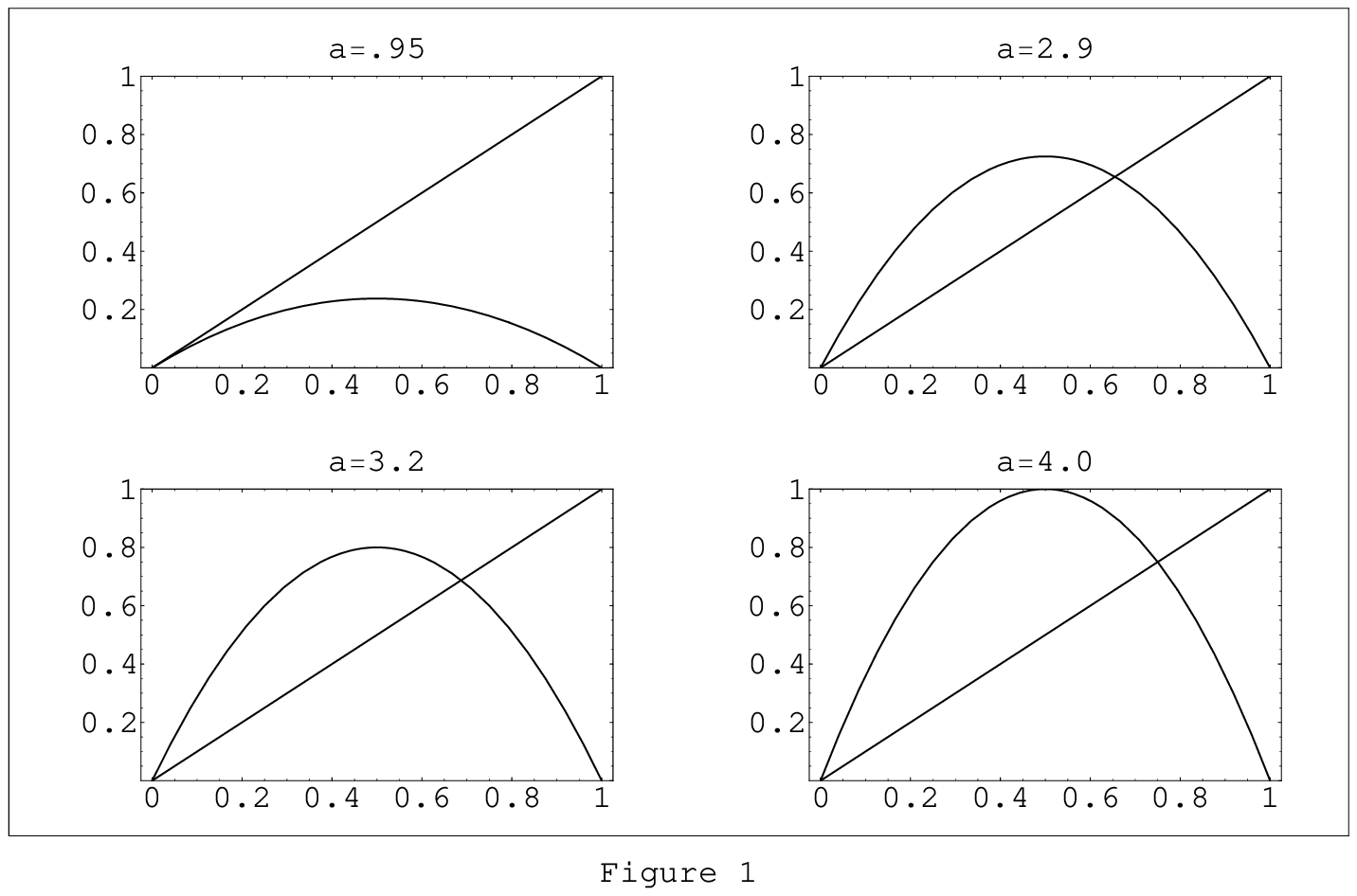} \caption{Graphs of $x_n$ versus $x_{n+1}$ for different
 values of $a$ }\label{dis1}
\end{figure}

The discrete logistic map described by
the single difference equation
 \begin{equation}x_{n+1}=a x_{n}(1-x_{n}) \label{roderick1}\end{equation}
as mentioned earlier, determines the  future value  of the
variable $x_{n+1}$ at time-step n+1 from the past value at
time-step n. The time evolution of $x_{n}$ generated by this
algebraic equation exhibits an extraordinary transformation from
order to chaos as the parameter a, which measures the strength of
nonlinearity is increased \cite{rod}. \\Although nonlinear
difference equations of this type
 have been studied extensively as simple models for turbulence in
 fluids, they also arise naturally in the study of evolution of
 biological populations. \\For the purpose of illustration we
 consider the population of gypsy moths in the northern United
 States, which exhibits wild and unpredictable fluctuations from
 year to year. However we could equally  well consider the
 evolution of economic prices determined by a nonlinear web
 model. \\Writing (\ref{roderick1}) in a slightly different form \[x_{n+1}=a x _{n}-a x^{2}_{n}\]
 we see that it is a simple quadratic equation, with the first
 term linear and the second term nonlinear. If the parameter
 $a>1$, the population increases, if $a<1$, the population
 decreases. If $a>1$, the population will eventually grow to a large
 enough value for the nonlinear term \ \ $-ax^{2}_{n}$,\ \ to become
 important. Since this term is negative, it represents a nonlinear
 death rate which dominates when the population is  too
 large. Biologically this nonlinear death rate could be due to the
 depletion of food supplies or the outbreak of  diseases in an
 over-crowded environment. The dynamics of this map and the
 dependence on the parameter a which measures the rate of linear
 growth and the size of the nonlinear term, are best understood
 using graphical analysis. Consider the graphs of $x_{n}$ versus
 $x_{n+1}$displayed in the Fig.(\ref{dis1}) for four different values of a.
 Equation (\ref{roderick1}) defines an inverted parabola with
 intercepts at \[x_{n}=0 \; and \; 1 \]  and a maximum value of \[x_{n+1}=
 a/4\] at \[x_{n}=0.5\]. \\Using these maps we can we can get a
 quantitative understanding of the dynamics of the logistic map
 in a quick way. \\Briefly, this graphical analysis
 tells us that if the normalized population starts out larger than
 1, then it immediately goes negative , becoming extinct in one
 time-step. Moreover if \ $a>4$, \  the peak of the parabola will exceed
 1, which makes it possible for initial populations near $0.5$ to
 become extinct in two time-steps. We will therefore restrict our
 analysis to values of  $a$  between 0 and 4. \\For values of $a<1$, the
 population always decreases to 0, as shown for $a=0.95$ in Fig.(\ref{dis1})
  The intersection of the parabola with the $45^{\tiny{o}}$ line
 at $x_{n}=0$ represents a stable fixed point on the map. Because
 $a$ is small a perturbation can be used to verify that almost all
 initial populations are attracted to this fixed point and become
 extinct. However for $a>1$ this fixed point becomes
 unstable. Instead the parabola now intersects the $45^{\tiny{o}}$
 line at \[x= \frac{a-1}{a}\]which corresponds to a new fixed
 point. \\For values  of \ $a$ \ between 1 and 3 almost all initial
 populations evolve to this equilibrium population. Then, as \ $a$ \ is
 increased between 3 and 4 , the dynamics change in remarkable
 ways. First the fixed point becomes unstable and the population
 evolves to a dynamic steady state in which it alternates between
 a large and small population. A time sequence converging to such a
 period-2 cycle is displayed in Fig.(\ref{dis1}) for a= 3.2 : the
 population cycles between two points on the parabola, $x_{n}\sim
 0.5$ and $x_{n}\sim 0.8$, in alternate years. For somewhat larger
 values of \ $a$\ , this period-2 cycle becomes unstable and is replaced
 by a period-4 cycle in which the population alternates
 high-low, returning to its original value every four time-steps. As
 \ $a$ \ is increased the longtime motion converges to
 period -- 8,16,32,64, cycles, finally accumulating to a cycle of
 infinite period $a=a_{\tiny{inf}}\sim 3.57$\\Having observed a period doubling
 sequence in numerical experiments Feigenbaum was able to prove
 that the intervals over which a  cycle is stable decreases at a
 geometric rate of $\sim 4.6692016$. The tremendous significance of
 this work is that this rate and other properties of the
 period-doubling bifurcation are universal in the sense  that they
 appear in the dynamics of any system which can be approximately
 modelled by a nonlinear map with a quadratic extremum.
Feigenbaum's theory has subsequently been confirmed by a wide
variety of physical systems such as turbulent fluids, oscillating
chemical reactions, nonlinear electrical circuits, and ring
lasers.
\\The investigation of period doubling in nonlinear dynamical
systems provides a superb example of the interplay between
numerical experiments and analytical theory. However, this
 sequence of regular periodic orbits is only the precursor to
 chaos. Included below is a bifurcation diagram , showing the
 beginning of Chaos:\begin{center}
\includegraphics{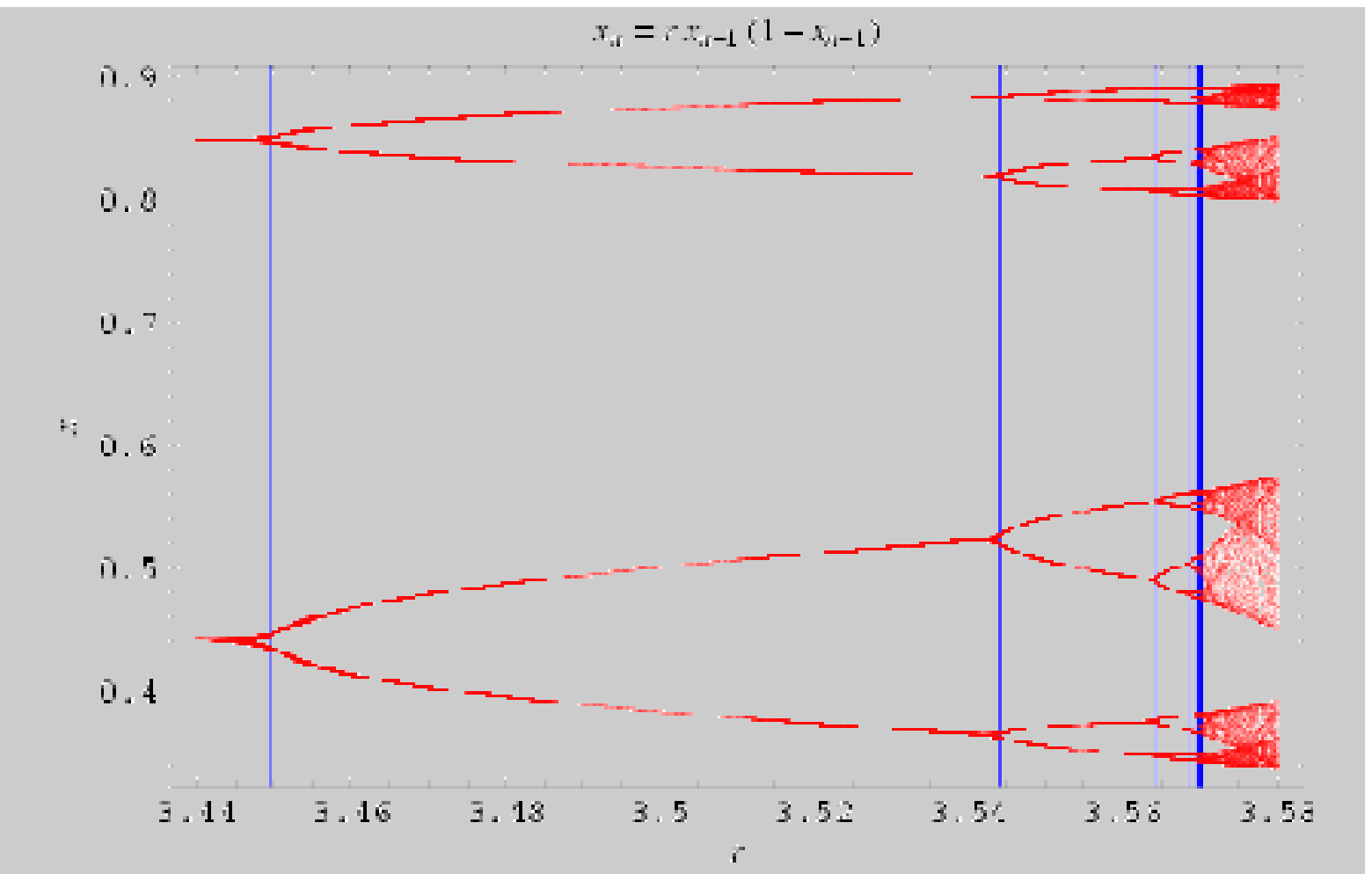}
\end{center}
\section{A Remarkable Mathematical Equivalence}
We would now like to deduce a mathematical equivalence between the
continuous and discrete cases and will comment on this. Let us
consider the equation from the continuous growth
\begin{equation}\label{diseq1}\frac{dP}{dt} = f(P), \hspace{2in} set \; \; P \equiv x
\end{equation} which gives
\begin{equation}\label{diseq2} x = \int f(x(t)) dt \end{equation}
where the integral is over suitable limits. Here f generalizes the
dependence of the right side of
\begin{equation}\label{logisticeq} \frac{d(P(t))}{dt}=r(M-P(t))P(t)
\hspace{2in} r > 0\end{equation} where P(t) is the population at
given time t and M is the maximum sustainable population
 \cite{r1,r2}
 on the population P at that time. Let us solve
equation (\ref{diseq2}) by the method of successive approximation
ie., we try on the right side of (\ref{diseq2}) a tentative
solution $x_{0}(t)$. (\ref{diseq2}) can then be written in the
form
\begin{equation}\label{diseq3} x_1 = F(x_0)\equiv \int f(x_0)dt \end{equation}
where $x_1$ gives the next level of approximation. Before
proceeding further, we remark that (\ref{diseq1}) is in the form
of the initial value problem,where the Lipschitz condition is
necessary for the convergence of the iterative procedure
\cite{coddington}. In a similar manner we get from (\ref{diseq3})
the more general equation
\begin{equation}\label{diseq4} x_{n+1} = F(x_{n}). \end{equation} (\ref{diseq4})
 can immediately
 identified with the discrete logistic
equation for example (\ref{roderick1}). It must be stressed
however that the discrete equation is based on a completely
different foundation that is the subscript $n$ in the discrete
case represents the population at the $n^{th}$ generation, whereas
in (\ref{diseq4}) $x_n$ represents the $n^{th}$ iteration or
approximation of the population $x \equiv P$ of the continuous
case. Nevertheless, this mathematical equivalence enables us to
apply the conclusions of the discrete case including the domain of
chaos. Thus chaotic behaviour of $x_n$ of the discrete case would
represent the lack of convergence of the iterates of the
continuous case. In this specific example if
\begin{equation}F(x_n) = a x_n (1-x_n),\label{diseq5} \end{equation}
then for $a$ = 3.57 the above iterative procedure breaks down.In
this case, using (\ref{diseq5}), it follows from (\ref{diseq3})
that
\[\int f(x) dt = a x (1-x).\]
Finally we remark that for a more conventional approach to the
above problem reference can be made to Krempasky, \cite{krem}.

\end{document}